\documentclass{dcds}
\usepackage{amsmath,amsthm,amssymb}
\usepackage{latexsym}
\usepackage{eucal}
  \usepackage[colorlinks=true]{hyperref}
\hypersetup{urlcolor=blue, citecolor=red}

\def\currentvolume{23}
 \def\currentissue{3}
  \def\currentyear{2009}
   \def\currentmonth{March}
    \def\ppages{1--XX}
     
\newtheorem{theorem}{Theorem}[section]
\newtheorem{lemma}{Lemma}[section]

\theoremstyle{definition}

\title[ Infinite superlinear growth of the gradient\ldots]
{Infinite superlinear growth of the gradient for the two-dimensional Euler equation}
\author[Sergey A. Denisov]{}
\keywords{Two-dimensional Euler equation, growth of the gradient}
 \subjclass{Primary 76B99, secondary 76F99.}
\email{denissov@math.wisc.edu}

\begin{document}
\maketitle

\centerline{\scshape Sergey A. Denisov} \medskip {\footnotesize
  \centerline{University of Wisconsin-Madison,
Mathematics Department}
  \centerline{480 Lincoln Dr. Madison, WI 53706-1388, USA}
}

\bigskip

 \centerline{(Communicated by Roger Temam)}

\begin{abstract}
For two-dimensional Euler equation on the torus, we prove that the
$L^\infty$ norm of the gradient can grow superlinearly for some
infinitely smooth initial data. We also show the exponential
growth of the gradient for finite time.
\end{abstract} \vspace{1cm}

\section{Introduction.\\}

 In this note, we are dealing with two-dimensional
Euler equation. We will write the equation for vorticity in the
following form

\[
\dot\theta=\nabla\theta \cdot u, \quad u=\nabla^\bot
\zeta=(\zeta_y,-\zeta_x), \quad \zeta=\Delta^{-1}\theta, \quad
\theta(x,y,0)=\theta_0(x,y)
\]
and $\theta$ is $2\pi$--periodic in both $x$ and $y$ (e.g., the
equation is considered on the torus $\mathbb{T}$). We assume that
$\theta_0$ has zero average over $\mathbb{T}$ and then
$\Delta^{-1}$ is well-defined since the Euler flow is
area-preserving and the average of $\theta(\cdot,t)$ is zero as
well.

 The global existence of the smooth solution for smooth
initial data is well-known \cite{March}. It is also known that the
gradient does not grow faster than the double exponential. On the
other hand, the lower bounds for the gradient's norm are not
studied well. There are some results on the infinite (not faster
than linear) growth of the gradient for a domain with the boundary
\cite{Yud,Yud1} or an annulus \cite{Nad}.

We will prove the following results. The proofs are inspired by
the recent preprint \cite{KN}.
\begin{theorem}
There is $\theta_0(x,y)\in C^\infty(\mathbb{T})$ such that
\begin{equation}
\frac{1}{T^2}\int\limits_0^T \|\nabla \theta(\cdot,t)\|_\infty
dt\to +\infty \label{superl}
\end{equation}
as $T\to+\infty$.\label{theorem1}
\end{theorem}
That implies the infinite superlinear growth of the time average
of the gradient.

\begin{theorem}
For any $T>0$, there is $\theta_0(x,y)\in C^\infty(\mathbb{T})$
with $\|\nabla \theta_0(x,y)\|_\infty\leq 10$, such that $\|\nabla
\theta(x,y,t)\|_\infty >0.1\exp(T/2)$ for some $t<T$.
\end{theorem}
That shows the possibility of the exponential growth over any
(fixed) period of time.

\section{Infinite superlinear growth.\\}

We start with several lemmas. We always assume that
$\theta_0(x,y)$ is infinitely smooth and its average over
$\mathbb{T}$ is zero. The following result is borrowed from the
recent preprint by A. Kiselev and F. Nazarov \cite{KN}.

\begin{lemma}
 Decompose
$\theta(x,y,t)=P_1 \theta+P_2\theta$, where $P_1$ is the Fourier
projector to the unit sphere on $\mathbb{Z}^2$.  If $\|P_2\theta
(\cdot,0)\|_2\leq \epsilon$, then $\|P_2\theta(\cdot,
t)\|_2\lesssim \epsilon$ for any $t>0$.\label{l1}
\end{lemma}
\begin{proof}
The following invariants are well-known
\[
\int_\mathbb{T} \theta^2dxdy=C_1, \int_\mathbb{T} \theta \zeta
dxdy=C_2
\]
Subtraction gives
\[
\sum\limits_{n_1^2+n_2^2>1}\left(1-\frac{1}{n_1^2+n_2^2}\right)|\hat\theta(n,t)|^2=const
\]
on the Fourier side. Since
\[
\left(1-\frac{1}{n_1^2+n_2^2}\right)\sim 1
\]
outside the unit sphere, we have the statement of the lemma.
\end{proof}
We also need the following elementary result on the preservation
of some symmetries.

\begin{lemma}
If $\theta_0$ is even (i.e. $\theta_0(x,y)=\theta_0(-x,-y)$, then
so is $\theta(\cdot, t)$. If $\theta_0$ is invariant with respect
to rotation by $\pi/2$ degrees around the origin, then so is
$\theta(\cdot,t)$.\label{sym1}
\end{lemma}
\begin{proof}
Assume that $\theta(x,y,t)$ is the solution. We need to show that
$\psi_1(x,y,t)=\theta(-x,-y,t)$ and $\psi_2(x,y,t)=\theta(-y,x,t)$
both are  solutions as well. Then, the uniqueness of the solution
to Cauchy problem would yield the statement of the lemma. Notice
that
\[
\dot\psi_1={\psi_1}_x\tilde{\zeta}_y-{\psi_1}_y\tilde{\zeta}_x
\]
where $\tilde{\zeta}(x,y)=\zeta(-x,-y)$. But  we also have
\[
\zeta(-x,-y)=\Delta^{-1}\psi_1
\]
as can be easily verified on the Fourier side. Therefore, $\psi_1$
solves Euler equation. For $\psi_2$,
\[
\dot \psi_2={\psi_2}_x\breve{\zeta}_y-{\psi_2}_y\breve{\zeta}_x
\]
where
\[
\breve{\zeta}(x,y,t)=\zeta(-y,x,t)
\]
On the other hand,
\[
\zeta(-y,x,t)=\sum_{n\in
\mathbb{Z}^2}e^{-in_1y+in_2x}\frac{\hat\theta(n_1,n_2,t)}{n_1^2+n_2^2}
\]
\begin{equation}
\Delta^{-1}\psi_2(x,y,t)=\sum_{n\in
\mathbb{Z}^2}e^{in_1x+in_2y}\frac{\hat{\psi_2}(n_1,n_2,t)}{n_1^2+n_2^2}\label{f2}
\end{equation}
\[
\hat\psi_2(n_1,n_2,t)=(2\pi)^{-2}\int_\mathbb{T}e^{-in_1x-in_2y}\theta(-y,x,t)dxdy
=\hat\theta(-n_2,n_1,t)
\]
after obvious change of variables. Changing indices $j_1=-n_2,
j_2=n_1$ in (\ref{f2}), we get $\breve{\zeta}=\Delta^{-1}\psi_2$
and therefore $\psi_2$ solves the Euler equation.
\end{proof}
Obviously, the symmetries considered follow immediately from the
symmetries of the multiplier $(n_1^2+n_2^2)^{-1}$. Another
symmetry preserved is $\theta(x,y,t)=-\theta(y,x,t)$. The proof of
this fact is similar but we are not going to use it. The following
elementary lemma will be used later
\begin{lemma}
If $a_j>0$ and
\[
\sum\limits_{j=1}^\infty a_j<\infty
\]
then
\[
\frac{1}{N^2}\sum\limits_{j=1}^N a_j^{-1}\to +\infty
\]
as $N\to+\infty$.\label{lemma3}
\end{lemma}
\begin{proof}We have
\[
\min\limits_{x_j> 0, x_1+\ldots+x_n\leq
\sigma}\sum\limits_{j=1}^nx_j^{-1}=\sigma^{-1}n^2
\]
Since
\[
\tau_N=\sum\limits_{j=N/2}^N a_j\to 0
\]
as $N\to\infty$, we have
\[
\frac{1}{N^2}\sum\limits_{j=1}^N
a_j^{-1}\geq\frac{1}{N^2}\sum\limits_{j=N/2}^N a_j^{-1}\geq
\frac{1}{4\tau_N}\to+\infty
\]
as $N\to +\infty$.
\end{proof}

Before giving the proof to the theorem 1, we notice that the
function $\theta^*(x,y,t)=\cos(x)+\cos(y)$ has spectrum on the
unit sphere. Therefore, $\zeta=\theta^*$ and $\theta^*$ is a
stationary solution. It is even and invariant with respect to
rotation by $\pi/2$ degrees. The flow generated by it is
\[
\nabla^\bot \zeta=(-\sin y,\sin x)
\]
and it is hyperbolic at the points $D_j=\pi(j_1,j_2)$ where
$j_{1}\mod 2\neq j_{2}\mod 2$. One can consider, e.g., the points
$A_1=(\pi,0)$, $A_2=(0,\pi)$ and their $2\pi$--translates. The
proofs will be based on certain stability of $\theta^*$ and
perturbation theory around the hyperbolic scenario which (without
nonlinear term) leads to infinite exponential growth of the
gradient for suitably chosen initial data. The idea is this: we
will show that away from the points $D_j$ the direction of the
flow is basically the same as without nonlinearity. Therefore, if
we place the bump at $D_j$, the area around $D_j$ and inside the
level set of this bump will be gradually exhausted by the flow.
But since the total area is preserved, this will be manifested
through narrowing of the ``chanel of exhaustion" thus leading to
collapse of two level sets and therefore growth of the gradient.

\noindent{\it Proof of theorem 1.} Let $\delta>0$ be small and
$U_\delta$ be the disc of radius $\sqrt\delta$ centered at origin.
Denote its $2\pi$--translates on $\mathbb{R}^2$ be $\hat{U}_\delta$.
Consider the new orthogonal coordinate system with the origin at
$O_1=A_1=(\pi,0)$ and axes $\xi: y=x-\pi$, $\eta: y=\pi-x$, the
orientation is positive. The relation to the original coordinates is
then

\[
\xi=\frac{y+x-\pi}{\sqrt 2}, \eta=\frac{y-x+\pi}{\sqrt 2}
\]
In this coordinate system, take the rectangle $\Pi_\delta=\{
|\xi|<0.1, |\eta|<\delta\}$. Rotate it around the origin by
$\pi/2$ degrees in the original coordinate system  and denote the
rectangle obtained by $\Pi'_\delta$. Consider all
$2\pi$--translates of $\Pi_\delta$ and $\Pi'_\delta$ and denote
the collection of all rectangles obtained this way by
$\hat{\Pi}_\delta$.

Consider $\theta_0(x,y)$ having the following properties:
\begin{itemize}
\item[(1)] $\theta_0$ is $2\pi$--periodic on $\mathbb{R}^2$, even, and is invariant
with respect to $\pi/2$--rotation around the origin.

\item[(2)] $\theta_0(x,y)=\theta^*(x,y)$ outside $\hat{U}_{\delta}$ and
$\hat\Pi_\delta$.

\item[(3)] In new coordinates $(\xi, \eta)$, $\theta_0=f(\xi,\eta)$ in
$\Pi_\delta$. The function $f$, defined on $\Pi_\delta$, is even,
$-1\leq f\leq 4$, the level set for $f=4$ is the segment
$\{\eta=0, |\xi|\leq 0.08\}$, the level set for $f=3$ is an
ellipse
\[
\left(\frac{\xi}{0.09}\right)^2+\left(\frac{\eta}{\delta/2}\right)^2=1
\]
\item[(4)] Inside $U_\delta$, we let $\theta_0=\theta^*-\phi_\delta$,
where $\phi_\delta\in C^\infty_0(U_\delta)$, is positive,
spherically symmetric, and such that
\[
\int_{\mathbb{T}}\theta_0(x,y)dxdy=0
\]
\item[(5)] $\theta_0(x,y)$ is infinitely smooth. That can be achieved by the obvious
smooth welding along the edge of $\Pi_\delta$.
\item[(6)] $\|\theta_0\|_\infty\leq 10$.

\end{itemize}
The existence of such $\theta_0$ is obvious. Then, the solution
$\theta(x,y,t)$ exists globally and is infinitely smooth. Due to
the lemma \ref{sym1}, it is even, $2\pi$--periodic, and is
invariant under the $\pi/2$--rotation around the origin.
Therefore, it must also be even with respect to all points $D_j$.
The function $\zeta(x,y,t)$ is therefore also even with respect to
the origin and $D_j$. If so, $\nabla^\bot \zeta=0$ at the origin
and at $D_j$, so we have $\theta(O,t)=const$,
$\theta(D_j,t)=const$, and these points do not move under the
flow.

 Write
$\theta=\theta^*+\psi$. By construction,
\[
\|\psi(\cdot,0)\|_2 \lesssim \delta^{1/2}
\]
and lemma \ref{l1} gives
\[
\|P_2\theta(\cdot,t)\|\lesssim \delta^{1/2}
\]

 Since $\|\theta^*\|_2=2\pi$, we have
\begin{equation}
|\|\theta(\cdot,0)\|_2-2\pi|\lesssim \delta^{1/2}\label{f3}
\end{equation}
and, because the $L^2$ norm is preserved by the flow, we have
(\ref{f3}) for any time. Therefore,
\begin{equation}
\sum\limits_{n_1^2+n_2^2=1}|\hat\theta(n_1,n_2,t)|^2=1+\underline{O}(\delta^{1/2})\label{f4}
\end{equation}
The symmetries for $\theta$ yield
\[
\hat\theta(n_1,n_2,t)=\hat\theta(-n_1,-n_2,t)=\hat\theta
(-n_2,n_1,t)
\]
for the Fourier coefficients.  That implies  $\hat\theta(n,t)$ is
a constant on the unit sphere and (\ref{f4}) shows that this
constant is $1/2+\underline{O}(\delta^{1/2})$. This, of course,
implies that
\[
\|\psi(\cdot,t)\|_2\lesssim \delta^{1/2}
\]

  Since $\|\theta_0\|_\infty\lesssim 1$ and all
$L^p(\mathbb{T})$ norms are preserved under the Euler flow, we
have that
\[
\|\psi(\cdot,t)\|_\infty\lesssim 1
\]
Interpolation gives
\begin{equation}
\|\psi(\cdot,t)\|_p\lesssim \delta^{\,1/p}, p>2 \label{f5}
\end{equation}
That essentially means $L^p$ stability of $\theta^*$, $p<\infty$.

Now, we need some simple perturbation estimates. We can write
\[
\nabla^\bot \zeta=(-\sin y,\sin x)+(g_1(x,y,t),g_2(x,y,t))
\]
where
\[
g_1=(\Delta^{-1} \psi)_y, g_2=-(\Delta^{-1}\psi)_x
\]
Since $\Delta^{-1}$ has the kernel with logarithmic singularity at
the origin, we get
\[
|g_{1(2)}(x,y,t)|\lesssim \int_{\mathbb{T}}
\frac{|\psi(\mu,\nu)|}{[(x-\mu)^2+(y-\nu)^2]^{1/2}}d\mu
d\nu\lesssim \frac{\|\psi\|_p}{(2-q)^{1/q}}\lesssim
\frac{\delta^{1/p}}{(2-q)^{1/q}}
\]
where $\quad q^{-1}+p^{-1}=1$. Fix any small $\epsilon>0$ and
arrange for $\delta$ such that
\[
\|g_{1(2)}(\cdot,t)\|_\infty<0.001\epsilon
\]
From the symmetry, we also know that $g_{1(2)}$ are odd with
respect to the origin and points $D_j$. The Euler equation can be
solved by the method of characteristics which yields the equations
for the flow
\begin{equation}
\left\{
\begin{array}{l}
\dot x=\sin y-g_1(x,y,t)\\
\dot y=-\sin x-g_2(x,y,t)
\end{array}\right.
\end{equation}
They can be rewritten in $(\xi,\eta)$ coordinates as
\begin{equation}
\left\{
\begin{array}{l}
\dot \xi=\sqrt 2\cos(\eta/\sqrt2)\sin (\xi/\sqrt2)-(g_2+g_1)/\sqrt2\\
\dot \eta=-\sqrt 2\cos(\xi/\sqrt2)\sin
(\eta/\sqrt2)+(g_1-g_2)/\sqrt2
\end{array}\right.
\end{equation}
If $\alpha=\xi/\sqrt2, \beta=\eta/\sqrt2$, then
\begin{equation}
\left\{
\begin{array}{l}
\dot \alpha=\cos\beta\sin \alpha+\mu_1\\
\dot \beta=-\cos\alpha\sin \beta+\mu_2
\end{array}\right.\label{flow3}
\end{equation}
and $\|\mu_{1(2)}\|_\infty\leq 0.01\epsilon$. We have
\begin{lemma}
Consider the Cauchy problem (\ref{flow3}) with initial data
$\alpha_0, \beta_0$ at time $t_0$. Then, if $|\alpha_0|\leq
3\epsilon$ and $|\beta_0|\leq 0.1$, then $|\beta(t_0+\tau)|<0.1$.
If $3\epsilon>|\alpha_0|>2\epsilon$ and $|\beta_0|<0.1$, then
$|\alpha(t_0+\tau)|>3\epsilon$. We take $\tau=1$.\label{lll}
\end{lemma}
\begin{proof}
The inequality
\[
|\dot\alpha|<|\alpha|+0.01\epsilon, |\alpha_0|\leq 3\epsilon
\]
implies
\[
|\alpha(t)|<4\epsilon e^\tau, t\in (t_0, t_0+\tau)
\]
By taking $\epsilon$ small, we get $|\beta(t)|<0.1$ for $t\in
(t_0,t_0+\tau)$ from the second equation. Now we just need to
notice that for $\alpha_0>2\epsilon$, $\alpha(t)$ grows on
$(t_0,t_0+\tau)$ and
\[
\dot \alpha\geq 0.9\alpha-0.01\epsilon, \quad
t\in [t_0,t_0+\tau]
\]
which implies
\[
\alpha(t_0+\tau)>\alpha_0e^{0.9\tau}-0.01\epsilon
\int\limits_0^\tau e^{0.9(\tau-s)}ds>3\epsilon
\]
\end{proof}
Let $E(t_1,t_2)$ denote the Euler flow from time $t_1$ to $t_2$.
$E(t_1,t_2)$ is an area-preserving diffeomorphism, symmetric with
respect to points $D_j$. It is important to mention that the level
curve $\theta=3$ around $O_1$ (which is an ellipse originally)
will always be homeomorphic to ellipse.

Consider the following sets. Take the set of points inside the
$\theta=3$ level curve at $t=0$ (an ellipse). Consider its
intersection with $B_{3\epsilon}=\{|\alpha|\leq 3\epsilon,
|\beta|<0.1\}$. Denote this set by $S_0$. We write $S_0=S_0^1\cup
S_0^2$ where $S_0^1=S_0\cap B_{2\epsilon}$. Take $E(0,1)S_0$. It
will be inside the tube $|\beta|<0.1$. Intersect it with
$B_{3\epsilon}$ and take the simply-connected component of this
set containing the point $O_1$ (which does not move under the
flow). Denote this set by $S_1$. It has the following properties:
\begin{itemize}
\item[(1)] The boundary of $S_1$ consists of the part of
$\theta(\cdot,1)=3$ level curve and parts of the vertical segments
$\alpha=\pm 3\epsilon$.

\item[(2)] The area
$|S_1|\leq |S_0^1|=|S_0|-|S_0^2|$. That simply follows from the
lemma \ref{lll} since $S_0^2$ will be carried away from
$B_{3\epsilon}$ by the flow.

\item[(3)] $S_1$ contains the part of the level curve $\theta(\cdot,1)=4$ which
is symmetric with respect to the origin and connects the following
points: the origin and $P_{\pm}^1$, where $\alpha$-coordinate of
$P_{\pm}^1$ is $\pm 3\epsilon$, respectively (i.e., they lie on
the left and right sides of $S_1$). Indeed, it follows from the
fact that $B_{2\epsilon}$--part of this level curve for $t=0$ will
stay inside the tube $|\beta|<0.1$ and its edges ($P_{\pm }^0$ )
will be carried away from $B_{3\epsilon}$ at time $1$.

\item[(4)] From the previous property we get that in the decomposition
$S_1=S_1^1\cup S_1^2$, the set $S_1^2$ is not empty and has a
positive area.

\end{itemize}
We then inductively define the sets  $S_n$ for all times $t_n=n$.
All properties given above will hold true. In particular,
\[
|S_{n+1}|\leq |S_n|-|S_n^2|
\]
which implies
\begin{equation}
\sum\limits_{n=0}^\infty |S_n^2|<\infty
\end{equation}
For each $n$, $S_n^2$ is symmetric with respect to $O_1$. Consider
its right part. The sides $\alpha=2\epsilon, \alpha=3\epsilon$
contain the points of the level curve: $\theta=4$. Therefore, we
have the trivial estimate
\[
|S_n^2|\gtrsim \epsilon \|\nabla \theta(\xi_n,n)\|_\infty^{-1}
\]
where $\xi_n$ is some point inside $S_n^2$. The application of
lemma \ref{lemma3} now yields
\[
\lim_{N\to\infty}\frac{1}{N^2} \sum_{n=0}^N \|\nabla
\theta(\cdot,n)\|_\infty=+\infty
\]
The obvious modification of this argument gives (\ref{superl}).
 $\Box$

\noindent\noindent{\bf Remark.} By time scaling, we can show that
the initial norm of the gradient can be taken as small as we like.
We believe, though, that there is an exponential growth of the
gradient in our scenario, not just superlinear. For instance, the
Euler evolution of the bump in the exterior hyperbolic flow in
$\mathbb{R}^2$ allows much better estimates than superlinear. In
that case, a simple multiscale argument allows to ``cut out" more
and more weight from the domain around zero. This is due to the fact
that the $L^p$ norm of the solution around zero will decrease
substantially in time. That we can not guarantee for the periodic
case.

Notice that the shear flow typically yields only linear growth of
the gradient. We never succeeded in applying our method to
perturbation of the shear flow (which is generated by another
stable stationary solution, say, $\theta^*(x,y)=\cos x$).

Analogous argument works for the family of equations where
$\zeta=\Delta^{-\gamma}\theta$ and $1>\gamma>1/2$ (although we do
not know the global existence of solution). It is an interesting
question to extend this proof to $\gamma=1/2$ (the so-called SQG).

\section{Exponential growth over finite time.\\}
In this section, we will prove theorem~2.\\ \noindent{\it Proof of
theorem 2.} We need the following elementary perturbation lemma
\begin{lemma}Consider the following system of equations
\begin{equation}\left\{
\begin{array}{l}
\dot \alpha=\alpha(1+f_1(\alpha,\beta,t))+f_2(\alpha, \beta,t)\beta\\
\dot \beta=-\beta(1+g_1(\alpha,\beta,t))+g_2(\alpha,\beta,t)\alpha
\end{array}\right.\label{flow}
\end{equation}
where $f_{1(2)},g_{1(2)}$ are $C^1$--smooth and
$\|f_{1(2)}\|_\infty, \|g_{1(2)}\|_\infty<\epsilon\ll 1$. Then, in
every neighborhood of the origin there is a pair
$(\alpha_0,\beta_0)$ such that the solution to the corresponding
Cauchy problem satisfies:
\[|\alpha(t)|\leq e^{-t/2}|\beta_0|, |\beta(t)|\leq
e^{-t/2}|\beta_0|
\]\label{pot}
\end{lemma}
\begin{proof}

Consider two sectors $S_1=\{\beta>2|\alpha|\}$ and
$S_2=\{\beta>|\alpha|\}$ and take any smooth curve
$\gamma_0=\gamma(0)$ without self-intersections connecting the
sides of $S_1$ and lying inside $S_1$ (see Picture 1). Let us
control the evolution of this curve $\gamma(t)$ under the flow
given by (\ref{flow}). Clearly, $\gamma(t)$ is smooth at any
moment $t$. Take some point $\alpha_0,\beta_0$ on $\gamma(0)$ and
consider its trajectory $\alpha(t),\beta(t)$. The second equation
easily implies that until this trajectory leaves the sector $S_2$,
we have
\begin{equation}
\beta_0e^{-2t} \leq \beta(t)\leq \beta_0e^{-t/2} \label{geom}
\end{equation}
Also, $\beta(t)$ decreases. Next, take the endpoint on the curve
$\gamma(0)$ which lies on the right side of $S_1$. Let it have
coordinates $(\alpha_0,2\alpha_0)$. The first equation of
(\ref{flow}) shows that
\[
\alpha(t)/2<\dot \alpha(t)<2\alpha(t)
\]
until the corresponding trajectory $(\alpha(t),\beta(t))$ is
inside the sector $\Omega_+=\{\alpha\leq \beta\leq 2\alpha,
\alpha>0\}$. Clearly, $\alpha(t)$ increases within this time
interval and we have
\[
\alpha_0e^{-2t} \leq \alpha(t)\leq \alpha_0e^{-t/2}
\]
Therefore, we can conclude that $(\alpha(t),\beta(t))\in \Omega_+$
for $t\in (0,t_0)$, where $t_0=(2\ln 2)/3$. Analogous inequalities
hold true for the other endpoint of $\gamma(0)$, the corresponding
trajectory will not leave $\Omega_-=\{-\alpha\leq \beta\leq
-2\alpha, \alpha<0\}$ as long as $t\in (0,t_0)$. It is also easy
to see that all of $\gamma(t)$ will be inside $S_2$ for this time
interval. Therefore, we can take time $t=t_0$ and consider the
part of the curve $\gamma(t_0)$ which has no self-intersections,
connects the opposite sides of $S_1$, and lies inside $S_1$. This
is possible by simple topological argument since the endpoints of
$\gamma(t_0)$ are in the sectors $\Omega_{\pm}$. Let us call this
new curve $\gamma_1$. Then, we consider the evolution of
$\gamma_1$ repeating the same construction again and again. We
will obtain the sequence of curves $\gamma_n$ which are all inside
$S_1$. Obviously, $\gamma_n$ is a part of $\gamma(nt_0)$.
  Now, one can
easily construct the solution with needed properties by the
standard approximation argument. Indeed, for time $t=nt_0$
consider a point on the curve $\gamma_n$ with $\alpha=0$ (there
might be many of those). Solve the equation backward obtaining the
trajectory. Consider the functions $\alpha_n(t),\beta_n(t)$ given
by this trajectory up to $nt_0$ and let
$\alpha_n(t)=0,\beta_n(t)=\beta_n(nt_0)e^{-(t-nt_0)/2}$ for
$t>nt_0$ (see Picture 2). The functions $\alpha_n(t),\beta_n(t)$,
considered on $[0,\infty)$, have uniformly bounded
$H^1(\mathbb{R}^+)$ norms. Indeed, by construction and
(\ref{geom}),
\[
0<\beta_n(t)\leq \beta_n(0)e^{-t/2}, |\alpha_n(t)|\leq
\beta_n(0)e^{-t/2}
\]
where the last inequality follows from the fact that the
trajectory is inside the sector $S_2$. The estimates on the
derivative now easily follow from (\ref{flow}). By the Alaoglu
theorem, there are $\alpha(t),\beta(t)\in H^1(\mathbb{R}^+)$, such
that
\[
\alpha_{n_k}(t)\to \alpha(t), \beta_{n_k}\to \beta(t)
\]
weakly in $H^1(\mathbb{R}^+)$. The Sobolev embedding is compact
and so the weak convergence in $H^1[0,b]$ implies the uniform
convergence on $[0,b]$ for any $b>0$. We have
\[
\alpha_{n_k}(t)=\alpha_{n_k}(0)+\int\limits_0^t
\alpha_{n_k}(s)[1+f_1(\alpha_{n_k}(s),\beta_{n_k}(s),s)]ds
\]
\[
+\int\limits_0^t
\beta_{n_k}(s)f_2(\alpha_{n_k}(s),\beta_{n_k}(s),s)ds
\]

\[
\beta_{n_k}(t)=\beta_{n_k}(0)-\int\limits_0^t
\beta_{n_k}(s)[1+g_1(\alpha_{n_k}(s),\beta_{n_k}(s),s)]ds
\]
\[
+\int\limits_0^t
\alpha_{n_k}(s)g_2(\alpha_{n_k}(s),\beta_{n_k}(s),s)ds
\]
where $t\in [0,b]$ and $k$ is large. Taking $k\to\infty$, we see
that $\alpha(t),\beta(t)$ satisfy the integral equations and
therefore are solutions to (\ref{flow}). Obviously, we  also have
\[
0<\beta(t)\leq \beta(0)e^{-t/2}, |\alpha(t)|\leq \beta(0)e^{-t/2}
\]
\ifx\JPicScale\undefined\def\JPicScale{1}\fi \unitlength
\JPicScale mm
\begin{picture}(50,50)(0,0)
\linethickness{0.1mm} \put(25,0){\line(0,1){50}}
\linethickness{0.1mm} \put(25,0){\line(1,0){25}}
\linethickness{0.1mm} \put(0,0){\line(1,0){25}}
\linethickness{0.1mm}
\multiput(25,0)(0.12,0.12){208}{\line(1,0){0.12}}
\linethickness{0.1mm}
\multiput(0,25)(0.12,-0.12){208}{\line(1,0){0.12}}
\linethickness{0.1mm}
\multiput(25,0)(0.12,0.24){167}{\line(0,1){0.24}}
\linethickness{0.1mm}
\multiput(5,40)(0.12,-0.24){167}{\line(0,-1){0.24}}
\linethickness{0.3mm} \qbezier(10,30)(12.6,32.61)(14.41,33.81)
\qbezier(14.41,33.81)(16.21,35.02)(17.5,35)
\qbezier(17.5,35)(18.8,35.01)(20,34.41)
\qbezier(20,34.41)(21.2,33.8)(22.5,32.5)
\qbezier(22.5,32.5)(23.8,31.2)(25,30.59)
\qbezier(25,30.59)(26.2,29.99)(27.5,30)
\qbezier(27.5,30)(28.8,30.02)(30,28.81)
\qbezier(30,28.81)(31.2,27.61)(32.5,25)
\qbezier(32.5,25)(33.8,22.39)(34.41,21.19)
\qbezier(34.41,21.19)(35.01,19.98)(35,20)
\put(45,25){\makebox(0,0)[cc]{$\Omega_+$}}

\put(5,25){\makebox(0,0)[cc]{$\Omega_-$}}

\put(30,40){\makebox(0,0)[cc]{$S_1$}}
\put(47,42){\makebox(0,0)[cc]{$S_2$}}
\put(15,37){\makebox(0,0)[cc]{$\gamma_n$}}

\put(15,-4){\makebox(0,0)[cc]{Picture 1}}

\put(50,-2){\makebox(0,0)[cc]{$\alpha$}}
\put(27,50){\makebox(0,0)[cc]{$\beta$}}
\end{picture}
\ifx\JPicScale\undefined\def\JPicScale{1}\fi \unitlength
\JPicScale mm
\begin{picture}(55,50)(0,0)
\linethickness{0.1mm} \put(35,0){\line(1,0){30}}
\linethickness{0.1mm} \put(10,0){\line(1,0){25}}
\linethickness{0.1mm} \put(35,0){\line(0,1){50}}
\linethickness{0.1mm}
\multiput(35,0)(0.12,0.12){208}{\line(1,0){0.12}}
\linethickness{0.1mm}
\multiput(10,25)(0.12,-0.12){208}{\line(1,0){0.12}}
\linethickness{0.1mm}
\multiput(35,0)(0.12,0.24){167}{\line(0,1){0.24}}
\linethickness{0.1mm}
\multiput(15,40)(0.12,-0.24){167}{\line(0,-1){0.24}}
\linethickness{0.3mm} \qbezier(35,15)(32.38,17.59)(32.38,20)
\qbezier(32.38,20)(32.38,22.41)(35,25)
\qbezier(35,25)(37.63,27.59)(37.03,30)
\qbezier(37.03,30)(36.43,32.41)(32.5,35)
\qbezier(32.5,35)(28.59,37.61)(26.78,38.81)
\qbezier(26.78,38.81)(24.98,40.02)(25,40)

\put(35,0){\line(0,1){15}}

\put(23,43){\makebox(0,0)[cc]{$(\alpha_k(t),\beta_k(t))$}}
\put(20,-4){\makebox(0,0)[cc]{Picture 2}}

\put(64,-2){\makebox(0,0)[cc]{$\alpha$}}
\put(38,50){\makebox(0,0)[cc]{$\beta$}}
\end{picture}

\end{proof}
\noindent{\bf Remark.} The lemma \ref{pot} is local, i.e. the
functions $f_{1(2)}, g_{1(2)}$ need to be defined and smooth only
around the origin.

Take any large $T$ and consider the initial value
\[
\theta(x,y, 0)=\theta^*(x,y)+\phi_\epsilon(x,y)
\]
where $\phi_\epsilon$ is supported around points $D_j$ and
$2\pi$-- translates of $O$. Around each point $D_j$ it is given by
(in $\alpha,\beta$--coordinates)
\[
\epsilon\phi(\alpha\epsilon^{-1},\beta\epsilon^{-1})
\]
where $\phi$-- nonnegative spherically symmetric bump with
$\phi(0,0)=1$ and support inside the unit disc. Around the origin,
$\phi_\epsilon$ is a similar bump chosen such that
\[
\int_\mathbb{T} \theta(x,y,0)dxdy=0
\]
Clearly, we can arrange for
\[
\|\nabla \theta(\cdot ,0)\|_\infty<10
\]
We will chose $\epsilon(T)$ later. The solution $\theta$ will
always exist, will be even, invariant under $\pi/2$--rotation, and
\[
\theta(\cdot,t)=\theta^*(\cdot)+\psi(\cdot,t)
\]
where
\[
\|\psi(\cdot,t)\|_2\lesssim \epsilon^2
\]
as follows form lemma \ref{l1} and symmetries of the solution.
Assume now that the statement of the theorem is wrong and
$\|\nabla \psi(\cdot,t)\|_\infty<0.1 \exp(T/2)+2$ for all $t\in
[0,T]$. Since $\psi$ is even with respect to all points $D_j$, its
gradient there is zero and we can write
\begin{equation}
(\Delta^{-1}\psi)_\alpha(\alpha,\beta)=
\nabla(\Delta^{-1}\psi)_\alpha(\alpha',\beta')\cdot(\alpha,\beta)
\label{appr}
\end{equation}
The analogous formula holds for the derivative in $\beta$. Let us
estimate the second derivatives of $\Delta^{-1}\psi$. We consider,
say, $(\Delta^{-1}\psi)_{\alpha\beta}$, the others can be treated
similarly. Since $\Delta^{-1}$ has a kernel with $\ln |z_1-z_2|$
singularity, we get
\[
(\Delta^{-1}\psi)_{\alpha\beta}(\alpha,\beta)\sim
\int\limits_{(\alpha-\xi)^2+(\beta-\eta)^2<10}
\frac{(\alpha-\xi)(\beta-\eta)}{((\alpha-\xi)^2+(\beta-\eta)^2)^2}\psi(\xi,\eta)d\xi
d\eta=
\]
\[
=\int_{10>(\alpha-\xi)^2+(\beta-\eta)^2>\rho^2}
\frac{(\alpha-\xi)(\beta-\eta)}{((\alpha-\xi)^2+(\beta-\eta)^2)^2}\psi(\xi,\eta)d\xi
d\eta\quad+
\]
\[
\int\limits_{(\alpha-\xi)^2+(\beta-\eta)^2<\rho^2}
\frac{(\alpha-\xi)(\beta-\eta)}{((\alpha-\xi)^2+(\beta-\eta)^2)^2}
\left[\psi(\alpha,\beta)+\nabla\psi(\xi',\eta')\cdot
(\xi-\alpha,\eta-\beta)\right] d\xi d\eta
\]
The first term is not larger than
\[
C\rho^{-1}\|\psi\|_2\lesssim \rho^{-1}\epsilon^2
\]
in absolute value. By our assumption, the second term is dominated
by $Ce^{T/2}\rho$. Thus, all second derivatives of
$\Delta^{-1}\psi$ are bounded by $C(\rho^{-1}\epsilon^2+\rho
e^{T/2})$. By choosing $\rho=0.001C^{-1}e^{-T/2},
\epsilon=0.001\sqrt{\rho C^{-1}}$, we obtain
\[
\|D^2\Delta^{-1}\psi\|_\infty\leq 0.01
\]
for $t\in [0,T]$. Then, the representation (\ref{appr}) allows to
write equations for the flow (\ref{flow3}) as
\begin{equation}
\left\{
\begin{array}{l}
\dot \alpha=\cos\beta\sin \alpha+\alpha f_1+\beta f_2\\
\dot \beta=-\cos\alpha\sin \beta+\alpha g_1+\beta g_2
\end{array}\right.\label{flow1}
\end{equation}
where $f_{1(2)}(\alpha,\beta,t), g_{1(2)}(\alpha,\beta,t)$ are
uniformly smaller than $0.01$. The simple modification of the
argument from lemma \ref{pot} shows existence of the flow
trajectory that starts at the level set $\theta=\epsilon \gamma,
\gamma<1/2$ and approaches the origin exponentially  fast. We
therefore have $\|\nabla \psi\|_\infty >0.5e^{t/2}$ thus giving a
contradiction at $t= T$.
 $\Box$

\noindent{\bf Remark.} Analogous argument shows the infinite
exponential growth for the problems where
$\zeta=\Delta^{-\gamma}\theta$ and $\gamma>1$. For that case, it is
also easy to prove that the gradient can not grow faster than the
exponential, so the exponential growth is sharp.

{\bf Acknowledgements.} This research was supported by NSF Grant
DMS-0500177 and by Alfred P. Sloan Research Fellowship. I thank A.
Kiselev and F. Nazarov for introducing me to the subject and for
many useful discussions. I am also grateful to A. Kiselev whose
comments and remarks helped a lot in improving the original
version of the manuscript.

\medskip

Received February 2008; revised July 2008.

\medskip

 \end{document}